\documentclass{amsart}

\usepackage{amsmath}
\usepackage{amsfonts}
\usepackage{amssymb}
\usepackage{amscd}
\usepackage{amsthm}
\usepackage{enumitem}
\usepackage{amsrefs}
\usepackage{hyperref}

\theoremstyle{plain}
\newtheorem{theorem}{Theorem}[section]
\newtheorem{lemma}[theorem]{Lemma}

\theoremstyle{definition}

\newcommand{\R}{\mathbb{R}}

\newcommand{\N}{\mathbb{N}}
\newcommand{\Z}{\mathbb{Z}}
\newcommand{\pP}{\mathbb{P}}

\newcommand{\pS}{\mathbb{S}}

\newcommand{\sS}{\ensuremath{\mathcal{S}}}

\newcommand{\reg}{\ensuremath{\mathcal{R}}}

\newcommand{\func}[3]{#1\colon#2\rightarrow #3}

\newcommand{\sC}{\mathcal{C}}
\newcommand{\smooth}{\mathcal{C}^{\infty}}

\title{Approximation of maps into spheres by regulous maps}
\author{Maciej Zieli{\'n}ski}
\address{Institute of Mathematics, Faculty of Mathematics and Computer Science, Jagiellonian University\\
ul.  \L{}ojasiewcza 6, 30-348\\
Krak{\'o}w, Poland}
\email{Maciej.Zielinski@im.uj.edu.pl}

\begin{document}

\begin{abstract}
Let $X$ be a compact real algebraic set of dimension $n$. We prove that every Euclidean continuous map from $X$ into the unit $n$-sphere can be approximated by regulous map. This strengthens and generalizes previously known results.
\end{abstract}
\maketitle
\section{Introduction}

A recent direction of research in real algebraic geometry is to study intermediate classes of maps between continuous and regular maps. Such classes as continuous rational maps, stratified-regular maps and regulous maps (which often coincide) have been studied in a series of papers \cites{BKV, FR, FMQ, KN, KKK, KUCH, KUCH2, KUCH3, KUCH4, KUCH5, KUCH6, KUCH7, KK, KK2, KK3, KK4, KZ, MON, N, ZIE}. The aim of this note is to strengthen a certain related result of \cite{KK} and \cite{KUCH3}.

We begin by fixing some terminology. A \textit{real algebraic variety} is a locally ringed space isomorphic to some algebraic subset of $\R^n$, for some positive integer $n$, endowed with the Zariski topology and the sheaf of regular functions. It is worth recalling that this class is identical with the class of quasi-projective real varieties, for more detail and information see \cite{BCR}. A morphism of real algebraic varieties is called a \textit{regular map}. We will be also interested in the Euclidean topology of such varieties and this is the topology we will mean, unless explicitly stated otherwise, when using topological notions. By a smooth map we understand a map of class $\smooth$.

Let $X$ be a real algebraic variety. A \textit{stratification} $\sS$ of $X$ is a finite collection of pairwise disjoint Zariski locally closed subvarieties whose union is equal to $X$. A map $\func{f}{X}{Y}$ of real algebraic varieties is said to be \textit{regulous} if it is continuous and if there exists some stratification $\sS$ of $X$ such that $f\vert_S$ is a regular map for every $S\in\sS$. We denote the set of all regulous maps between $X$ and $Y$ by $\reg^0(X,Y)$. We shall treat $\reg^0(X,Y)$ as a subspace of the space $\sC(X,Y)$ of all continuous maps endowed with the compact-open topology. Note that regulous maps in the sense of our definition were called stratified-regular in \cite{KK} and the followup papers \cites{KK2, KK3}. This definition is different but equivalent to that of \cite{FR} where the terminology was introduced, see \cite{KK}*{Remark 2.3} or \cite{KZ}.

Each regulous map is also \textit{continuous rational} - i.e. $f$ is continuous and $f\vert_{X^0}$ is regular for some Zariski open dense subset $X^0\subset X$. While the converse is false in general, it is true if $X$ is nonsingular, see \cite{KN}.

Let us first recall the following result contained in \cite{KUCH3}.
\begin{theorem}\label{w1}
Let $X$ be a compact nonsingular real algebraic variety of dimension $p\geq 1$. Then the set $\reg^0(X,\pS^p)$ is dense in $\sC(X,\pS^p)$.
\end{theorem}
A related weaker result allowing for a singular $X$ is contained in \cite{KK}:
\begin{theorem}\label{w2}
Let $X$ be a compact real algebraic variety of dimension $p\geq 1$. Then any continuous map $X\rightarrow \pS^p$ is homotopic to a regulous map.
\end{theorem}
Our aim is to strenghten Theorem \ref{w2} by showing that the nonsingularity assumption of Theorem \ref{w1} is unnecessary.
\begin{theorem}\label{main}
Let $X$ be a compact real algebraic variety of dimension $p\geq 1$. Then the set $\reg^0(X,\pS^p)$ is dense in $\sC(X,\pS^p)$.
\end{theorem}

It is well-known that analogous results do not hold if regulous maps in Theorems \ref{w1}, \ref{w2}, and \ref{main} are  replaced with regular maps. For example, a continuous map $\pS^1\times\pS^1\rightarrow\pS^2$ is homotopic to a regular map if and only if it is null-homotopic, cf. \cite{BK}.

\section{Proof of the main theorem}

We shall use the concept of the algebraic cohomology classes of a real algebraic variety which we recall now. Let $X$ be a compact nonsingular real algebraic variety. A class in $H^*(X;\Z/2)$ is said to be algebraic if it is Poincar{\'e} dual to a homology class in $H_*(X,\Z/2)$ represented by an algebraic subset. The set $H^*_{\text{alg}}(X; \Z/2)$ of all algebraic cohomology classes is a subring of $H^*(X;\Z/2)$ and if $\func{f}{X}{Y}$ is a regular map, then the induced map $f^*$ in cohomology maps $H^*_{\text{alg}}(Y; \Z/2)$ into $H^*_{\text{alg}}(X, \Z/2)$, cf \cites{BT, BCR, BH}.

An important tool we need is \cite{KUCH3}*{Lemma 2.2}, which allows for controlled approximation of continuous maps into projective space by regular maps. We restate it here for convenience.

\begin{lemma}\label{kohom}
Let $X$ be a compact nonsingular real algebraic variety and let $A$ be a Zariski closed subvariety of $X$. Let $\func{f}{X}{\pP^n(\R)}$ be a continuous map whose restriction $\func{f\vert_A}{A}{\pP^n(\R)}$ is a regular map. Assume that \[f^*(H^1(\pP^n(\R);\Z/2))\subset H^1_{\mathrm{alg}}(X;\Z/2).\] Then one can find a regular $\func{g}{X}{\pP^n(\R)}$ arbitrarily close to $f$ and satisfying $g\vert_A=f\vert_A$ (i.e every neigborhood of $f$ in $\sC(X,\pP^n(\R))$ contains such a map).
\end{lemma}

We are now ready to prove Theorem \ref{main}

\begin{proof}
Let $f$ be any map in $\sC(X,\pS^p)$. Treating $X$ as a closed subset of $\R^m$ for some $m\in\N$, one can find a smooth map $\func{f_0}{U}{\pS^p}$ defined on some neighborhood $U$ of $X$ in $\R^m$ such that $f_0\vert_X$ is arbitrarily close to $f$. Let $\Sigma$ denote the singular locus of $X$. Then, $f_0(\Sigma)\subsetneq \pS^p$, since $\dim\Sigma<p$. This allows us to approximate $f_0\vert_\Sigma$ by regular maps using the stereographic projection and Weierstrass approximation theorem. We can therefore reduce the proof (by replacing $f$ with suitably modified $f_0$) to the case where $f$ is a restriction of a smooth map defined on a neighborhood of $X$ in $\R^m$ and $f\vert_\Sigma$ is regular with $f(\Sigma)\subsetneq\pS^p$. Then, by Sard's theorem, there exists an $s_0\in\pS^p\setminus f(\Sigma)$ which is a regular value of the smooth map $f\vert_{X\setminus \Sigma}$. 

By Hironaka's resolution of singularities theorem \cites{HIR, KOL}, there exists a finite composition of blowups $\func{\pi}{Y}{X}$ over $\Sigma$ with $Y$ nonsingular. The restriction $\func{\bar{\pi}}{Y\setminus\pi^{-1}(\Sigma)}{X\setminus\Sigma}$ of $\pi$ is then a biregular isomorphism and $s_0$ is a regular value of the smooth map $f\circ\bar{\pi}$. Letting $F=(f\circ\pi)^{-1}(s_0)=(f\circ\bar{\pi})^{-1}(s_0)$, consider the blowup of $Y$ with center $F$, which we will denote $\func{\sigma}{B(Y,F)}{Y}$, and the blowup $\func{\tau}{B(\pS^p, s_0)}{\pS^p}$ of $\pS^p$ over $s_0$. Since $\dim X = p$, the set $F$ is finite as a fiber of the smooth map $f\circ\bar{\pi}$ over its regular value $s_0$, hence $B(Y,F)$ is a real algebraic variety. Finiteness of $F$ and the fact that $f\circ\pi$ is also a smooth map allow us to apply \cite{AK}*{Lemma 2.5.9} in order to construct a smooth lifting $g$ of $f\circ\pi$ to a map between $B(Y,F)$ and $B(\pS^p, s_0)$ making the following diagram commute:
\[\begin{CD}
B(Y,F) @>g>> B(\pS^p, s_0)\\
@V\sigma VV     @V\tau VV \\
Y @>f\circ\pi >> \pS^p \\
@V\pi VV         @| \\
X @>f >> \pS^p
\end{CD}\]

Our aim for now is to find a regular map $\func{H}{B(Y,F)}{B(\pS^p,s_0)}$ arbitrarily close to $g$ in $\sC(B(Y,F), B(\pS^p,s_0))$ in such a way that the map $\func{\tilde{f}}{X}{\pS^p}$ making the following diagram commute will be regulous and close to $f$:
\[\begin{CD}
B(Y,F) @>H>> B(\pS^p, s_0)\\
@V\pi\circ\sigma VV     @V\tau VV \\
X @>\tilde{f}>>     \pS^p
\end{CD}\] 
Let $D=\sigma^{-1}(F)$ and $E=\tau^{-1}(s_0)$. Then $D$ and $E$ are nonsingular algebraic hypersurfaces in the respective blowups. Let $u\in H^1(B(Y,W);\Z/2)$ be the cohomology class Poincar{\'e} dual to the homology class in $H_*(B(Y,F);\Z/2)$ represented by $D$ and let $v$ be the class in $H^1(B(\pS^p,s_0);\Z/2)$ dual to the class in $H_*(B(\pS^p,s_0);\Z/2)$ represented by $E$. Recall that there exists a biregular isomorphism $\func{\varphi}{B(\pS^p, s_0)}{\pP^p(\R)}$ such that
\begin{equation}\label{projective}
\varphi(E)=\pP^{p-1}(\R)\subset \pP^p(\R).
\end{equation}
Hence $H^1(B(\pS^p,s_0);\Z/2)=\varphi^*(H^1(\pP^p(\R);\Z/2))\cong \Z/2$ with $v$ the generator. Since $g$ is transverse to $E$ and $D=g^{-1}(E)$ we have $u=g^*(v)$ (cf. \cite{BH}*{Proposition 2.15}) and it follows that 
\begin{equation}\label{g*}
g^*(H^1(B(\pS^p, s_0);\Z/2))\subset H^1_{\text{alg}}(B(Y;F);\Z/2).
\end{equation}

Let $\func{i}{D}{B(Y,F)}$ and $\func{j}{E}{B(\pS^p, s_0)}$ be the inclusion maps and let $\func{\bar{g}}{D}{E}$ be the map induced by the restriction of $g$ to $D$. Then, since $g\circ i=j\circ\bar{g}$ we have $\bar{g}^*(j^*(v))=i^*(g^*(v))=i^*(u)$. This allows us to apply Lemma \ref{kohom} to approximate $\bar{g}$ by regular maps. Indeed, since by (\ref{projective}) $H^1(E;\Z/2)$ is generated by $j^*(v)$ and by definition $i^*(u)$ is in $H^1_\text{alg}(D;\Z/2)$, we have
\[\bar{g}^*(H^1(E;\Z/2))\subset H^1_\text{alg}(D;\Z/2).\]
Therefore, there exists a regular map $\func{r}{D}{E}$ arbitrarily close to $\bar{g}$ in $\sC(D,E)$. 

Let $L=(\pi\circ\sigma)^{-1}(\Sigma)$ and let $\func{h}{\Sigma}{B(\pS^p,s_0)}$ be the regular map given by $(\tau\vert_{\tau^{-1}(\pS^p\setminus\{s_0\})})^{-1}\circ(f\vert_\Sigma)$. We now extend $r$ to a smooth map $\func{G}{B(Y,F)}{B(\pS^p, s_0)}$ close to $g$ such that $G\vert_L=(h\circ\pi\circ\sigma)\vert_L$ (this is possible since the sets $E$ and $L$ are disjoint). If $G$ is close enough to $g$ we have $g^*=G^*$ and hence by (\ref{g*}) we can apply Lemma \ref{kohom} to $G$ with $A=D\cup L$. This gives a regular map $\func{H}{B(Y,F)}{B(\pS^p, s_0)}$ which is close to $g$ and satisfies $H\vert_L=(h\circ\pi\circ\sigma)\vert_L$ and $H(D)\subset E$. The latter shows that the map $\func{\bar{H}}{Y}{\pS^p}$ given by:
$$\bar{H}(y)=\begin{cases}
s_0\text{, for }y\in F,\\
\tau(H(\sigma^{-1}(y)))\text{, for }y\notin F,\\ \end{cases}$$
is continuous and well-defined. Moreover $\bar{H}(y)=f\circ\pi(y)$ for $y\in\pi^{-1}(\Sigma)$. Also note that $\bar{H}$ is close to $f\circ\pi$, since $(\tau\circ H)\vert_{Y\setminus F}$ is close to $\tau\circ f\circ\pi$. Therefore, the map $\func{\tilde{f}}{X}{\pS^p}$ given by
$$\tilde{f}(x)=\begin{cases}
f(x)\text{, for }x\in \Sigma,\\
\bar{H}(\pi^{-1}(x))\text{, for }x\notin \Sigma,\\
\end{cases}$$
is well-defined, continuous and can be chosen arbitrarily close to $f$.

It remains to check that $\tilde{f}$ is regulous. Since $\tilde{f}\vert_{\Sigma}$ is regular as is $\tilde{f}$ in restriction to the finite set $f^{-1}(s_0)$, it is enough to show $\tilde{f}\vert_{X\setminus A}$ is regular where $A=\Sigma\cup f^{-1}(s_0)$. Since we have $\tilde{f}\vert_{X\setminus A}=\tau\circ H\circ(\sigma\vert_{Y\setminus F})^{-1}\circ(\pi\vert_{X\setminus \Sigma})^{-1}$ with all the maps on the right-hand side regular, the proof is finished. 
\end{proof}

\begin{subsection}*{Acknowledgements}
The author was partially supported by the National Science Centre (Poland) under grant number 2014/15/B/ST1/00046.
\end{subsection}


\begin{bibdiv}
\begin{biblist}
\bib{AK}{book}{
	author={Akbulut,S.},
	author={King, H.},
	title={Topology of Real Algebraic Sets},
	series={Math. Sci. Res. Inst. Publ.},
	volume={25},
	publisher={Springer},
	year={1992},
}


\bib{BT}{article}{
	author={Benedetti, R.},	
	author={Tognoli, A.},	
	title={Remarks and counterexamples in the theory of real algebraic vector bundles and cycles},
	book={
	title={G\'{e}om\'{e}trie alg\'{e}brique r\'{e}elle et formes quadratiques},
    series={Lecture Notes in math.},
	volume={959},
    publisher={Springer},
    year={1982},
    },
    pages={198--211},
    
}

\bib{BKV}{article}{
	author={Bilski, M.},	
	author={Kucharz, W.},
	author={Valette, A.},		
	author={Valette, G.},
	title={Vector bundles and regulous maps,},
	journal={Math. Z.},
    volume={275},
    year={2013},
    pages={403--418},

}

\bib{BCR}{book}{
	title={Real Algebraic Geometry},
	author={Bochnak, J.},
	author={Coste, M.},
	author={Roy, M.-F.},
	series={Ergeb. der Math. und ihrer Grenzgeb. Folge (3)},
	volume={36},
    publisher={Springer-Verlag, Berlin},
    year={1998},
}

\bib{BK}{article}{
	author={Bochnak, J.},
	author={Kucharz, W.},
	title={Realization of homotopy classes by algebraic mappings},
	journal={J.Reine Angew. Math.},
    volume={377},
    year={1987},
    pages={159--169},
}

\bib{BH}{article}{
	author={Borel, A.},
	author={Haefliger, A.},
	title={La classe d'homologie fondamentale d'ub espace analytique},
	journal={Bull. Soc. Math. France},
    volume={89},
    year={1961},
    pages={461--513},
}

\bib{FR}{article}{
	author={Fichou, G.},
	author={Huisman, J.},
	author={Mangolte, F.},
	author={Monnier, J.-Ph.},
	title={Fonctions r\'egulues},
	journal={J. Reine Angew. Math.},
    volume={718},
    year={2016},
    pages={103--151},
}
\bib{FMQ}{article}{
	author={Fichou, G.},
	author={Monnier, J.-P.},
	author={Quarez, R.},
	title={Continuous functions in the plane regular after one blowing up},
	journal={Math. Z.},
	doi = {10.1007/s00209-016-1708-8},
}

\bib{HIR}{article}{
	author={Hironaka, H.},
	title={Resolution of singularities of an algebraic variety over a field of characteristic zero},
	journal={Ann. of Math.},
	volume={79},
	number={2},
	pages={109--326},
	year={1964},
}

\bib{KOL}{book}{
	  author={Koll{\'a}r, J.},
      title={Lectures on resolution of singularities},
      series={Ann. of Math. Stud.},
      volume={166},
      year={2007},
      publisher={Princeton University Press}
}
    

\bib{KN}{article}{
	  author={Koll{\'a}r, J.},
	  author={Nowak, K.},
      title={Continuous rational functions on real and {$p$}-adic varieties},
      journal={Math. Z.},
      volume={279},
      year={2015},
      number={1-2},
      pages={85--97},
}

\bib{KKK}{article}{
	author={Koll{\'a}r, J.},
	author={Kucharz, W.},
	author={Kurdyka, K.},
	title={Curve-rational functions},
	journal={Math. Ann.},
	doi={10.1007/s00208-016-1513-z},
	year={2017},
}

\bib{KUCH}{article}{
	author={Kucharz, W.},
	title={Rational maps in real algebraic geometry},
	journal={Adv. Geom.},
	volume={9},
	year={2009},
	pages={517--539},
}

\bib{KUCH2}{article}{
	author={Kucharz, W.},
	title={Regular versus continuous rational maps},
	journal={Topology Appl.},
	volume={160},
	year={2013},
	pages={1086--1090},
}


\bib{KUCH3}{article}{
	author={Kucharz, W.},
	title={Approximation by continuous rational maps into spheres},
	journal={J. Eur. Math. Soc.},
	volume={16},
	year={2014},
	pages={1555--1569},
}

\bib{KUCH4}{article}{
	author={Kucharz, W.},
	title={Continuous rational maps into the unit 2-sphere},
	journal={Arch. Math. (Basel)},
	volume={102},
	year={2014},
	pages={257--261},
}

\bib{KUCH5}{article}{
	author={Kucharz, W.},
	title={Some conjectures on stratified-algebraic vector bundles},
	journal={J. Singul.},
	volume={12},
	year={2015},
	pages={92--104},
}

\bib{KUCH6}{article}{
	author={Kucharz, W.},
	title={Continuous rational maps into spheres},
	journal={Math. Z.},
	volume={283},
	year={2016},
	pages={1201--1215},
}

\bib{KUCH7}{article}{
	author={Kucharz, W.},
	title={Stratified-algebraic vector bundles of small rank},
	journal={Arch. Math. (Basel)},
	volume={107},
	year={2016},
	pages={239--249},
}

\bib{KK3}{article}{
	author={Kucharz, W.},
	author={Kurdyka, K.},
    title={Some conjectures on continuous rational maps into spheres},
    journal={Topology Appl.},
    volume={208},
    year={2016},
    pages={17--29},
}

\bib{KK}{article}{
	author={Kucharz, W.},
	author={Kurdyka, K.},
	title={Stratified-algebraic vector bundles},
	journal={J. Reine Angew. Math.},
	year={2016},
	doi={10.1515/crelle-2015-0105},
}

\bib{KK2}{article}{
	author={Kucharz, W.},
	author={Kurdyka, K.},
    title={Comparison of stratified-algebraic and topological K-theory},
    eprint={arXiv:1511.04238 [math.AG]},    
}

\bib{KK4}{article}{
	author={Kucharz, W.},
	author={Kurdyka, K.},
    title={Linear equations on real algebraic surfaces},
    eprint={arXiv:1602.01986 [math.AG]},    
}

\bib{KZ}{article}{
	author={Kucharz, W.},
	author={Zieli\'nski, M.},
    title={Regulous vector bundles},
    eprint={arXiv:1703.05566 [math.AG]},    
}

\bib{MON}{article}{
	author={Monnier, J.-P.},
    title={Semi-algebraic geometry with rational continuous functions},
    eprint={arXiv:1603.04193 [math.AG]},    
}

\bib{N}{article}{
	  author={Nowak, K.J.},
      title={Some results of algebraic geometry over Henselian rank one valued fields},
      journal={Sel. Math. New Ser.},
      year={2017},
      pages={455--495},
      volume={28},
}

\bib{ZIE}{article}{
	author={Zieli{\'n}ski, M.},
    title={Homotopy properties of some real algebraic maps},
    journal={Homology Homotopy Appl.},
    volume={18},
    year={2016},
    pages={287--294},    
}

\end{biblist}
\end{bibdiv}
\end{document}